\documentclass[12pt, leqno]{article}
\usepackage[english]{babel}
\usepackage[latin1]{inputenc}
\usepackage{amsmath}
\usepackage{amsfonts}
\usepackage{amssymb}
\usepackage{a4wide}
\newcommand{\Card}{{\rm Card}}

\newcommand{\C}{\mathbb{C}}

\newcommand{\Id}{{\rm Id}}

\newcommand{\K}{\mathbb{K}}

\newcommand{\N}{\mathbb{N}}

\newcommand{\Q}{\mathbb{Q}}
\newcommand{\Z}{\mathbb{Z}}

\newcommand{\al}{\alpha}

\newcommand{\congru}{\equiv}

\newcommand{\croix}{\times}

\newcommand{\eps}{\varepsilon}

\newcommand{\etoile}{^\ast}

\newcommand{\findem}{\end{pf}}

\newcommand{\moins}{\setminus}

\newcommand{\psd}{\rtimes}

\newcommand{\ro}{\varrho}

\newcommand{\La}{{\rm La}}

\newcommand{\zdzp}{(\Z / 2 \Z)^p}

\newcommand{\zdzd}{(\Z / 2 \Z)^2}
\newcommand{\zdz}{\Z / 2 \Z}

\newcommand{\spp}{{\mathfrak S}_p}

\newcommand{\sdeux}{{\mathfrak S}_2}

\newcommand{\jsoul}{\underline{j}}
\newcommand{\ssoul}{\underline{s}}

\newcommand{\epssoul}{\underline{\eps}}
\newcommand{\epsprsoul}{\underline{\eps'}}

\newcommand{\zeron}{\{0,\ldots,n\}}
\newcommand{\unp}{\{1, \ldots, p\}}
\newcommand{\unpmu}{\{1, \ldots, p-1\}}
\newcommand{\unq}{\{1, \ldots, q\}}
\newcommand{\zerop}{\{0, \ldots, p\}}

\newcommand{\sunp}{s_1, \ldots, s_p}
\newcommand{\sunq}{s_1, \ldots, s_q}
\newcommand{\sunj}{s_1, \ldots, s_j}
\newcommand{\sunk}{s_1, \ldots, s_k}
\newcommand{\runp}{r_1, \ldots, r_p}
\newcommand{\junp}{j_1, \ldots, j_p}
\newcommand{\xiunp}{\xi_1, \ldots, \xi_p}
\newcommand{\kunp}{k_1, \ldots, k_p}
\newcommand{\Kunp}{K_1, \ldots, K_p}
\newcommand{\ellunp}{\ell_1, \ldots, \ell_p}
\newcommand{\sigunp}{\sigma_1, \ldots, \sigma_p}

\newcommand{\epsunp}{\eps_1, \ldots, \eps_p}
\newcommand{\epsprunp}{\eps'_1, \ldots, \eps'_p}
\newcommand{\epsgammuunp}{\eps_{\gamma^{-1}(1)}, \ldots, \eps_{\gamma^{-1}(p)}}

\newtheorem{Th}{Theorem}[section]
\newtheorem{Prop}[Th]{Proposition}
\newtheorem{Cor}[Th]{Corollary}

\newtheorem{Conj}[Th]{Conjecture}
\newtheorem{Exemple}{\it Example\/}
\newtheorem{Remarque}{\it Remark}

\def\Dem{\hspace{-\parindent}{{\bf Proof }}}

\newcommand{\Ic}{I^{{\rm c}}}

\newcommand{\VK}{V_{\K}}
\newcommand{\VprK}{V'_{\K}}
\newcommand{\WK}{W_{\K}}
\newcommand{\WprK}{W'_{\K}}

\newcommand{\WprQ}{W'_{\Q}}

\newcommand{\Cvarpi}{C[\varpi]}
\newcommand{\Cprvarpi}{C'[\varpi]}
\newcommand{\Cvidesj}{C [ \emptyset, \emptyset, (s_1, \ldots, s_p) , (j_1, \ldots, j_p)]}
\newcommand{\Cprvidesj}{C' [ \emptyset, \emptyset, (s_1, \ldots, s_p) , (j_1, \ldots, j_p)]}

\newcommand{\Cprvidesjsoul}{C' [ \emptyset, \emptyset, \ssoul, \jsoul]}
\newcommand{\Cvidesjsun}{C [ \emptyset, \emptyset, (1,s_2,  \ldots, s_p) , (j_1,j_2,  \ldots, j_p)]}
\newcommand{\Cvidesjgam}{C [ \emptyset, \emptyset, (s_{\gamma(1)}, \ldots, s_{\gamma(p)}) , (\eps_1 \cdot j_{\gamma(1)}, \ldots, \eps_p \cdot j_{\gamma(p)}  )]}
\newcommand{\Cvidesprjpr}{C [ \emptyset, \emptyset,  (s'_1, \ldots, s'_p), (j'_1, \ldots, j'_p)]}

\newcommand{\Domsigz}{D[\varpi, \sigma_1, \ldots, \sigma_q; z]}

\newcommand{\zetaf}{\zeta_{\rm f}}
\newcommand{\zetaantisym}{\zeta^{{\rm as}}}  

\newcommand{\calz}{{\cal Z}}
\newcommand{\cala}{{\cal A}}

\newcommand{\tildero}{\tilde{\ro}}
\newcommand{\sigmati}{\tilde{\sigma}}

\title{Multiple series connected to Hoffman's conjecture
on multiple zeta values}
\author{S. Fischler\footnote{Univ. Paris-Sud, Laboratoire de Math\'ematiques, UMR CNRS 8628, 
B\^atiment 425, 91405 Orsay, France;  stephane.fischler@math.u-psud.fr}} 
\date{\today}

\begin{document}
\maketitle

\setcounter{tocdepth}{2}
\baselineskip 6mm

\begin{abstract}
Recent results of Zlobin and Cresson-Fischler-Rivoal allow one to decompose any suitable $p$-uple series of hypergeometric type into a linear combination (over the rationals) of multiple zeta values of depth at most $p$; in some cases, only the multiple zeta values with 2's and 3's are involved (as in Hoffman's conjecture). In this text, we study the depth $p$ part of this linear combination, namely the contribution of the multiple zeta values of depth exactly $p$. We prove that it satisfies some symmetry property as soon as the $p$-uple series does, and make some conjectures on the depth $p-1$ part of the linear combination when $p=3$. Our result generalizes the property that (very) well-poised univariate hypergeometric series involve only zeta values of a given parity, which is crucial in  the proof by Rivoal and Ball-Rivoal that $\zeta(2n+1)$ is irrational for infinitely many $n \geq 1$.
The main feature of the proof is an algebraic approach, based on representations of $\zdzp \psd \spp$.
\end{abstract}

\bigskip


\section{Introduction}

The {\em multiple zeta values} (also called {\em multiple harmonic series}) are defined, for integers $s_1 \geq 2$ and $s_2, \ldots, s_p \geq 1$ (with $p \geq 0$), by
$$\zeta(s_1, \ldots, s_p) = \sum_{k_1 > \ldots > k_p \geq 1} \frac{1}{k_1^{s_1} \ldots k_p ^{s_p}}.$$
The integer $p$ is called the {\em depth}, and $s_1 + \ldots + s_k$ the {\em weight}. 
They appear in many areas of mathematics, and are related to motives, knots, renormalization, \ldots . Many linear or algebraic  relations (over $\Q$) are known between these numbers (see for instance \cite{MiW} or \cite{ZudilinMZV} for a survey), so that it is not easy to find a basis for the vector space $\calz$ spanned over $\Q$ by the multiple zeta values. Hoffman has made the following conjecture \cite{Hoffman}, which seems to be completely out of reach nowadays:
%,  for which numerical evidence is known up to weight 16 ? :

\begin{Conj} \label{conjho} 
The multiple zeta values $\zeta(s_1, \ldots, s_p)$, with $p \geq 0$ and $s_i \in \{2, 3\}$ for any $i \in \unp$, make up a basis of $\calz$.
\end{Conj}

The following  (apparently weaker) conjecture is actually equivalent to  Hoffman's:

\begin{Conj} \label{conjli} 
The multiple zeta values $\zeta(s_1, \ldots, s_p)$, with $p \geq 0$ and $s_i \in \{2, 3\}$ for any $i \in \unp$, are linearly independent over $\Q$.
\end{Conj}

 Indeed, assume that Conjecture \ref{conjli} holds.  
For any $n \geq 0$, let $\calz_n$ be the  subspace of  $\calz$ spanned by the multiple zeta values of weight $n$, $\cala_n$ be the set of all sequences $(s_1, \ldots, s_p ) \in \{2,3\}^p$ with $p \geq 0$ and $s_1 + \ldots + s_p = n$, and $\cala = \cup_{n \geq 0} \cala_n$. For instance, $\calz_0 = \Q$ (since $\zeta(\emptyset) = 1$), $\calz_1 = \{0\}$, $\calz_2 = \Q \zeta(2) = \Q \pi^2$ (since $\zeta(2) = \pi^2/6$), $\calz_3 = \Q \zeta(3)$ (since $\zeta(2,1) = \zeta(3)$);
$\cala_0 = \{\emptyset\}$, $\cala_1 = \emptyset$, $\cala_2 = \{ (2) \}$, $\cala_3 = \{ (3) \}$. 

The upper bound $\dim_\Q \calz_n \leq \Card \,  \cala_n$ has been proved for any $n \geq 0$ (\cite{Terasoma}, see also \cite{Goncharov}). Now Conjecture \ref{conjli} implies the converse inequality, so that equality holds (as conjectured by Zagier), and   the $\zeta(s_1, \ldots, s_p)$, for $(s_1, \ldots, s_p) \in \cala_n$, make up a basis of $\calz_n$. Moreover,  Conjecture \ref{conjli} implies that $\calz$ is the direct sum of the $\calz_n$; this is enough to deduce Conjecture \ref{conjho}. 

\bigskip

The goal of this paper is to provide tools (following from an algebraic point of view) for proving partial results towards Conjecture \ref{conjli} or other linear independence conjectures about multiple zeta values. A usual way to prove that some numbers are linearly independent over $\Q$ is to produce very small (but non-zero) linear forms, with not too big integer coefficients, in these numbers. This allows one  to apply a linear independence criterion, for instance Nesterenko's \cite{Nest}. 

\bigskip

If we restrict our attention to depth 1, that is values $\zeta(s)$ of Riemann $\zeta$ function at integers $s \geq 2$, one conjectures that $1$, $\zeta(2)$, $\zeta(3)$, $\zeta(4)$, \ldots, are linearly independent over $\Q$. One could hope to prove this using the following fact, which can be proved easily using partial fraction expansion and in which $(k)_{\al} = k(k+1)\ldots(k+\al-1)$ denotes Pochhammer's symbol:

\medskip

$(i)$ For any integers $r, n \geq 0$ and $A \geq 2$, and any polynomial $P \in \Q[k]$ of degree at most $A(n+1)-2$, the real number 
\begin{equation} \label{eqre}
\sum_{k \geq 1}\frac{P(k)}{(k+r)_{n+1}^A}
\end{equation}
is a linear combination, with rational coefficients, of  $1$, $\zeta(2)$, $\zeta(3)$,  \ldots, $\zeta(A)$.

\medskip

Actually this hope is completely out of reach for the moment, and all one can hope for is partial results towards the conjecture. Since $\zeta(2k)$ is a rational multiple of $\pi^{2k}$ for any $k \geq 1$, the most interesting point concerns values $\zeta(s)$ for odd integers $s$. To obtain a diophantine result about these values, one may use the following result due to Rivoal \cite{RivoalCRAS} and Ball-Rivoal \cite{BR}:

\medskip

$(ii)$ With the notation of $(i)$, if 
\begin{equation} \label{eqsymBR}
P(-n-k) = (-1)^{A(n+1)+1} P(k)
\end{equation}
then \eqref{eqre} is a linear combination, with rational coefficients, of  $1$ and the  $\zeta(s)$ for odd values of $s$ with $3 \leq s \leq A$.

\medskip

Using appropriate choices of $P$ satisfying \eqref{eqsymBR}, Rivoal \cite{RivoalCRAS} and Ball-Rivoal \cite{BR} were able to deduce diophantine results from  Nesterenko's linear independence criterion \cite{Nest}, for instance:

\medskip

$(iii)$ The $\Q$-vector space spanned by $1$, $\zeta(3)$, $\zeta(5)$, $\zeta(7)$, \ldots is infinite-dimensional.

\bigskip

Now let us come back to multiple zeta values of arbitrary depth. The aim we have in mind is to prove diophantine results (in the style of $(iii)$) towards Conjecture \ref{conjli} (though our results are more general, and may be useful for other conjectures about multiple zeta values). The first step in this direction (analogous to $(i)$) is the following   theorem of Zlobin  \cite{ZlobinZametki2005}, which gives a large family of series which are linear forms,  with rational coefficients, in the multiple zeta values of the conjecture.

\begin{Th} \label{th1} 
Let $p \geq 1$ and $  n_1, \ldots, n_p,  r_1, \ldots, r_p$ be non-negative integers such that 
\begin{equation} \label{eqsepth1}
r_i \geq r_{i+1} + n_{i+1} + 1  \mbox{ for any } i \in \unpmu .
\end{equation}
Let   $P(k_1, \ldots, k_p)$ be any polynomial with rational coefficients such that $\deg_{k_i}P \leq 3n_i + 1$ for any $i \in \unp$.   Then the series
\begin{equation} \label{eqth1}
\sum_{k_1\geq \ldots \geq k_p \geq 1}
\frac{P(k_1,  \ldots, k_p)}{(k_1+r_1)_{n_1+1}^{3} (k_2+r_2)_{n_2+1}^{3}\ldots (k_p+r_p)_{n_p+1}^{3}}
\end{equation}
is a  linear combination,  with rational coefficients,  of  multiple zeta values $\zeta(s_1, \ldots, s_q)$ with
$q \in \zerop$ and $s_i \in \{2, 3\}$ for any $i \in \unq$. 
\end{Th}

However, Conjecture \ref{conjli} is still out of reach. For a given integer $p \geq 1$, we are also very far from knowing how to prove that $\delta_p = 2^{p+1}-1$, where $\delta_p$ is the dimension of the vector space spanned (over $\Q$) by the multiple zeta values $\zeta(s_1, \ldots, s_q)$ with
$q \in \zerop$ and $s_i \in \{2, 3\}$ for any $i \in \unq$. Since $\zeta(2,2,\ldots, 2)=\pi^{2q}/(2q+1)!$, the transcendence of $\pi^2$ yields $\delta_p \geq p+1$. Any improvement in this lower bound seems out of reach for the moment, so it might be useful to consider vector spaces spanned only by specific such multiple zeta values (see however \cite{Zlobinpoidsimpair}). Proving lower bounds for the dimension of such $\Q$-vector spaces seems to be the only reasonable hope for obtaining new results towards  Conjecture \ref{conjli}. For instance the linear independence of 1, $\zeta(2)$, and $\zeta(3)$ is still an open question (eventhough Ap\'ery \cite{Apery} has proved in 1978 that $\zeta(3)$ is irrational). The following weaker statements already seem to be very difficult conjectures:

\begin{Conj} \label{conj13} 
Among the numbers $1$, $\zeta(2)$, $\zeta(3)$,  and $\zeta(2,3) - \zeta(3,2)$, at least three are linearly independent over $\Q$. 
\end{Conj}

\begin{Conj} \label{conj14} 
Among the numbers 1, $\zeta(2)$, $\zeta(3)$, $\zeta(2,3)$, $\zeta(3,2) $, $\zeta(3,3) $, and  $\zeta(3,2,3)$,  at least three are linearly independent over $\Q$. 
\end{Conj}

\bigskip

These conjectures could play the role of assertion $(iii)$ above. In this paper, we make a step towards  Conjecture \ref{conj13}  by proving the corresponding assertion $(ii)$ (see Theorem \ref{th3} below). We also state a conjectural assertion $(ii)$ corresponding to Conjecture \ref{conj14} (see Conjecture \ref{conjintro} below). Of course Conjecture \ref{conj14} is weaker than Conjecture \ref{conj13}; the crucial point in Conjecture \ref{conj14} is that $\zeta(2,2)$ is not involved (since $1$, $\zeta(2) = \pi^2/6$ and $\zeta(4) = \pi^4/90$ are known to be linearly independent).

\bigskip

More generally, in the present text we state and prove multivariate statements analogous to assertion $(ii)$, namely we  refine a generalization of  Theorem \ref{th1} so that  (under suitable assumptions) only a restricted set of multiple zeta values appear in the linear combination. We mostly achieve this goal with respect to the depth $p$ part of this  linear combination. For instance, we prove the following result:

\begin{Th}  \label{th2} 
In the setting of Theorem \ref{th1}, let $\sigunp \in \{2,3\}$ and assume 
that
 $$P(k_1, \ldots, k_{i-1}, -k_i -2 r_i - n_i , k_{i+1}, \ldots, k_p) = (-1)^{n_i+1+\sigma_i} P(\kunp)$$
 for any $i \in \unp$. Then in the linear combination that represents \eqref{eqth1}, we may assume that $\zeta(\sigunp)$ is the only multiple zeta value of depth $p$ that appears with a (possibly) non-zero coefficient.
 \end{Th}

 In fact our main result involves not only symmetry properties $k_i \mapsto  -k_i -2 r_i - n_i $ (as in Ball-Rivoal's statement $(ii)$), but also permutations of the variables $k_1 + r_1 + \frac{n_1}{2}$, \ldots, $k_p + r_p + \frac{n_p}{2}$. To state it, we define    two actions of the group $G = \zdzp \psd \spp$: one on rational functions 
 \begin{equation} \label{eqRint}
 R(\kunp) = \frac{P(k_1,  \ldots, k_p)}{(k_1+r_1)_{n_1+1}^{A_1}  \ldots (k_p+r_p)_{n_p+1}^{A_p}} ,
 \end{equation}
 and the other one on symbols $\zetaf(\sunp)$ corresponding to multiple zeta values (but with no linear relations between them). 
 This algebraic approach is the main feature of this text. Our main result (see \S \ref{subsecenplus}) reads as follows: if a subgroup $H$ of $G$ acts on  $R(\kunp)$ through a character $\chi$, then it acts  in the same way on  the depth $p$ part of the linear combination of multiple zeta values that represents
\begin{equation} \label{eqintro}
\sum_{k_1 \geq \ldots \geq k_p \geq 1}R(\kunp).
\end{equation}
 Here and throughout this text, the {\em depth $k$ part} of such a linear combination 
$$\sum_{j=0} ^p \sum_{s_1, \ldots, s_j} \lambda[\sunj] \zeta(\sunj)$$
is 
$$\sum_{s_1, \ldots, s_k} \lambda[\sunk] \zeta(\sunk). $$
The relation with Hoffman's conjecture appears only when $A_1 = \ldots = A_p = 3$, but we treat the general situation (so that assertion $(ii)$ above is a special case of our results, in the easy case $p=1$). The subgroup $H = \zdzp \croix \Id$ is used to deduce Theorem \ref{th2}. When $p=2$ and $H = \{(1,1)\} \croix \sdeux$, the following statement can  be obtained (as a special case of  Theorem \ref{thcyclic} below):

\begin{Th}  \label{th3} 
Let $n, r_1, r_2 \geq 0$ be integers such that $r_1 \geq r_2 + n + 1$. Let  $P \in\Q[k_1,k_2]$ be any polynomial, of degree at most $3n+1$ with respect to each variable, such that
\begin{equation} \label{eqPsymint}
P(k_2 + r_2 - r_1 + \frac{n_2-n_1}{2}, k_1 + r_1 - r_2 + \frac{n_1-n_2}{2}) = -P(k_1,k_2).
\end{equation}
Then the series 
\begin{equation} \label{eqPsymintso}
\sum_{k_1\geq k_2  \geq 1}
\frac{P(k_1,k_2)}{(k_1+r_1)_{n+1}^{3} (k_2+r_2)_{n+1}^{3}}
\end{equation}
is a linear combination over the rationals of $1$, $\zeta(2)$, $\zeta(3)$ and $\zeta(2,3) - \zeta(3,2)$. 
 \end{Th}

This theorem could be a tool to prove Conjecture \ref{conj13}. It would be sufficient to construct polynomials $P$ such that \eqref{eqPsymintso} is very small but non-zero, and the coefficients of the 
 linear combination are not too big. Then Nesterenko's linear independence criterion \cite{Nest} would give the result. However, we have no idea of appropriate choices for $P$.
 
 \bigskip
 
Our results enable us to construct linear combinations of multiple zeta values of depth at most $p$, with a good control upon the depth $p$ part. In some cases, we can even make it vanish (see Remark \ref{rem35} in \S \ref{subsecdep2}). But the   main drawback is  that we describe only the depth $p$ part of the linear combination arising from a $p$-uple series. However, we are confident that they can be extended in some way to the depth $p-1$ part (and maybe further ?), at least when $p=3$. This has been done in a   special case in \cite{CFRsym} (see Example \ref{exemplesym} in \S \ref{subsecenplus}), but the proof is very complicated and does not use any algebraic structure. An interesting challenge would be to prove the result of   \cite{CFRsym} with the same kind of algebraic methods as the ones introduced here.

As far as extensions to the  depth $p-1$ part when $p=3$ are concerned,  we have checked the following conjecture (using the algorithm \cite{CFRweb}) for $n \leq 2$.
 
 \begin{Conj} \label{conjintro}
Let $n$ be a non-negative integer. Denote by $\sigma$ the element of $\{2,3\}$ such that $\sigma \congru n \bmod 2$, and by $\sigmati$ the other element of $\{2,3\}$. Let $P(k_1, k_2, k_3)$ be a polynomial, with rational coefficients, such that
$$
\left\{
\begin{array}{l}
P(-k_1 -5n-4 , k_2, k_3) = -P(k_1, k_2, k_3)\\
P(k_1, - k_2 - 3n-2, k_3) = +P(k_1, k_2, k_3)\\
P(k_1,  k_2 ,  - k_3- n) = -P(k_1, k_2, k_3)
\end{array}
\right.
$$
and $\deg_{k_i} P \leq 3n+1$ for any $i \in \{1,2,3\}$. 
Then the  series
$$\sum_{k_1 \geq k_2 \geq k_3 \geq 1} 
\frac{ P(k_1, k_2, k_3) }{(k_1+2n+2)_{n+1}^3 (k_2+n+1)_{n+1}^3 (k_3)_{n+1}^3}
$$
 is a linear combination (over the rationals) of  1, $\zeta(2)$, $\zeta(3)$, $\zeta(2,3)$, $\zeta(3,2) $, $\zeta(3,3) $, and  $\zeta(\sigma, \sigmati, \sigma)$. 
\end{Conj}
In this assertion, the depth 1 part of the linear combination follows from Theorem \ref{th1} and the depth 3 part from Theorem \ref{th2}, so that the only conjectural aspect is the shape of the depth 2 part, namely the fact that $\zeta(2,2)$ does not appear in it.

\bigskip

Proving this conjecture and finding good choices for $P$ could lead to a proof of Conjecture \ref{conj14},  in the same way as Theorem \ref{th3} could be used to prove Conjecture \ref{conj13}. 

\bigskip

The structure of this text is as follows. We state in Section \ref{secprelim} the known results about the expansion of a multiple series \eqref{eqintro} as a linear form in multiple zeta values,
in a general setting of which Theorem \ref{th1} is only a special case. We also refine these results and state  consequences on the so-called ``derivation procedure''. Section \ref{secprinc} is devoted to the statement and proof of our main result, together with corollaries obtained in special cases. At last, Section \ref{secconj} deals with conjectures and open questions.

\bigskip

\noindent {\bf Acknowledgements: } The statements we prove or conjecture here have been guessed thanks to many computations using GP/Pari on the Medicis computing centre. 

\section{General Results} \label{secprelim}

In this section, we first summarize in Theorem \ref{thzl} (\S \ref{subsecresuenonces}) the known results about multiple series of hypergeometric type, proved in \cite{ZlobinZametki2005} and   \cite{CFRalgo}. Then we prove (\S \ref{subsecraff}) a refinement concerning the depth $p$ part in this theorem, and an easy consequence of it (namely a generalization of the derivation procedure used classically in depth 1, see \S \ref{subsecderiv}). 

\subsection{Decomposition of a Multiple Series} \label{subsecresuenonces}

Let $p \geq 1$ and $A_1, \ldots, A_p, n_1, \ldots, n_p,  r_1, \ldots, r_p$ be non-negative integers.
Let   $P(k_1, \ldots, k_p)$ be any polynomial. We let for $j \in \{1, \ldots, p\}$:
$$
D_j = \Big( \sum_{i=1} ^j A_i (n_i+1) \Big) - j - 1, 
$$
and we assume 
 \begin{equation} \label{eqCV}
\sum_{i=1}^j \deg_{k_i}P \leq D_j \mbox{ for any } j \in \{1, \ldots, p\}.
\end{equation}
This condition is necessary and sufficient for the series \eqref{eq7} below to converge
 (see \cite{CFRalgo}, \S 8.4). 

\medskip

Let $\K$ be a subfield of $\C$. The important case is $\K = \Q$, since the general case follows from it by linearity (considering monomials $P$). Other fields $\K$ will be used in \S \ref{subseccyclic}.

\begin{Th} \label{thzl}
Let $p \geq 1$ and $A_1, \ldots, A_p, n_1, \ldots, n_p,  r_1, \ldots, r_p$ be non-negative integers.
Let   $P(k_1, \ldots, k_p)$ be any polynomial with  coefficients in $\K$,  such that the conditions \eqref{eqCV} hold. 
Then the series
\begin{equation} \label{eq7}
\sum_{k_1\geq \ldots \geq k_p \geq 1}
\frac{P(k_1,  \ldots, k_p)}{(k_1+r_1)_{n_1+1}^{A_1} (k_2+r_2)_{n_2+1}^{A_2}\ldots (k_p+r_p)_{n_p+1}^{A_p}}
\end{equation}
is a linear combination (over $\K$) of  multiple zeta values
$\zeta (s_1 ,\dots ,s_q )$ with $0\leq q\leq p$ and $\sum_{j=1}^q s_j \leq \sum_{j=1}^p A_j$.
More precisely, we may assume that if $\zeta (s_1 ,\dots ,s_q )$ appears with a non-zero coefficient then there exist $1 = \ell_0 < \ell_1 < \ldots < \ell_q = p+1$ such that $s_j \leq A_{\ell_{j-1}} + A_{\ell_{j-1}+1} + \ldots + A_{\ell_j-1}$ for any $j \in \unq$. 

If we assume also
\begin{equation} \label{eq6}
r_i \geq r_{i+1} + n_{i+1} + 1  \mbox{ for any } i \in \unpmu 
\end{equation}
then \eqref{eq7} is a  linear combination (over $\K$)  of  multiple zeta values $\zeta(s_1, \ldots, s_q)$ with:
\begin{itemize}
\item  $0 \leq q \leq p$, 
\item  $s_1 \geq 2$,
\item there exist $1 \leq i_1 < \ldots < i_q \leq p$ with $1 \leq s_{\ell}  \leq A_{i_\ell}$ for any $\ell \in \unq$,
\item $\Card \,   \{\ell \in \unq, \, s_\ell = 1\} \leq \Card \,  \{ i \in \unp, \, \deg_{k_i} P \geq A_i(n_i+1)-1\}$. 
\end{itemize}

As a corollary, if \eqref{eq6} holds and 
\begin{equation} \label{eq8}
\deg_{k_i} P \leq A_i(n_i+1)-2 \mbox{ for any } i \in \unp
\end{equation}
then \eqref{eq7} is a  linear combination (over $\K$) of   $\zeta(s_1, \ldots, s_q)$ with $0 \leq q \leq p$  for which there exist $1 \leq i_1 < \ldots < i_q \leq p$ with $2 \leq s_{\ell}  \leq A_{i_\ell}$ for any $\ell \in \unq$.  
\end{Th}

The first part of this  theorem was proved independently in  \cite{ZlobinZametki2005} (Theorem 1) and \cite{CFRalgo} (Th\'eor\`eme 3); the second part (in which \eqref{eq6} is assumed to hold)  follows from the proof of \cite{ZlobinZametki2005}  (Theorem 5).

\begin{Remarque} \label{remunicite}
Since there are many linear relations over $\Q$ among multiple zeta values, the linear combination in Theorem \ref{thzl} is not unique. However, the proof produces a specific one; so does\footnote{But it is not clear to us whether both ways produce the same  linear combination.} the algorithm \cite{CFRweb}. Throughout this text, when we claim that ``the'' linear combination in Theorem \ref{thzl}  satisfies some additional property, we refer to the one constructed from the proof. The important point, in general, is merely the existence of such a linear combination with the additional property.
\end{Remarque}

It should be possible to prove a general statement, assuming that $r_i \geq r_{i+1} + n_{i+1} + 1 $ for any $i$ in a subset $I$ of $\unpmu$, such that the cases $I = \emptyset $ and $ I = \unpmu$ give the two parts of Theorem \ref{thzl}. 

\bigskip

In the statement of Theorem \ref{thzl}, the multiple zeta value 1 appears as $\zeta(s_1, \ldots, s_q)$ with $q=0$. 

When \eqref{eq6} and \eqref{eq8} hold, $\K = \Q$, and $A_i = 3$ for any $i$, Theorem \ref{thzl} reduces to Theorem \ref{th1} stated in the Introduction. In the case where  $\K = \Q$, \eqref{eq6}  holds, $r_i  = r_{i+1} + n_{i+1} + 1$ and neither $n_i$ nor $A_i$ depend on $i$, we obtain for instance the following Corollary.

\begin{Cor} \label{cornonagglo}
Let $p \geq 1$, $A \geq 0$, $n \geq 0$ be integers. Let $P(k_1, \ldots, k_p)$ be any polynomial with rational coefficients such that the conditions \eqref{eqCV} hold. Then the series
\begin{equation} \label{eqcornonagglo}
\sum_{k_1\geq \ldots \geq k_p \geq 1}
\frac{P(k_1,  \ldots, k_p)}{(k_1+(p-1)(n+1))_{n+1}^{A} (k_2+(p-2)(n+1))_{n+1}^{A}\ldots (k_p)_{n+1}^{A}}
\end{equation}
is a  linear combination,  with rational coefficients,  of  multiple zeta values $\zeta(s_1, \ldots, s_q)$ with $0 \leq q \leq p$, $s_1 \geq 2$, $1 \leq s_i \leq A$ for any $i \in \unq$, such that the number of $i$ with $s_i =1$ is not greater than the number of $i$ with $\deg_{k_i}P \geq A(n+1)-1$. 

Accordingly, if  
$\deg_{k_i} P \leq A(n+1)-2 $ for any $ i \in \unp$
then \eqref{eqcornonagglo} is a  linear combination,  with rational coefficients,  of   $\zeta(s_1, \ldots, s_q)$ with $0 \leq q \leq p$  and $2 \leq s_i \leq A$ for any $i \in \unq$. 
\end{Cor}

\subsection{A Refinement of the Depth $p$ Part} \label{subsecraff}

 In one variable, for integers $e, A, n, r \geq 0$, we have the partial fraction expansion
$$\frac{k^e}{(k+r)_{n+1}^A} = \sum_{f=0}^{e-A(n+1)}D_f k^f + \sum_{j=0}^n \sum_{s=1}^A \frac{C_{j,s}}{(k+r+j)^s}$$
with rational numbers $D_f$ and $C_{j,s}$ (the first sum does not appear if $e-A(n+1)<0$). Applying this identity with respect to $k_1$, \ldots, $k_p$, we obtain the following partial fraction expansion:
\begin{equation} \label{eqvarpi}
\frac{P(k_1,  \ldots, k_p)}{(k_1+r_1)_{n_1+1}^{A_1}  \ldots (k_p+r_p)_{n_p+1}^{A_p}} = 
\sum_{\varpi} 
 \Cvarpi \frac{\prod_{i \in I} k_i^{f_i}}{\prod_{i \in \Ic}(k_i+r_i+ j_i)^{s_i}}.
\end{equation}
In this formula, we denote by $J$ the set of all indices $i$ such that $\deg_{k_i} P  \geq A_i(n_i+1)$; $I$ is a subset of $J$, $\Ic = \unp \moins I$, $(f_i)_{i \in I}$ is a family of non-negative integers such that $f_i \leq \deg_{k_i} P - A_i(n_i+1)$ for any $i \in I$, and $(s_i)_{i \in \Ic}$ and $(j_i)_{i \in \Ic}$ are families of non-negative integers such that $1 \leq s_i \leq A_i$ and $0 \leq j_i \leq n_i$  for any $i \in \Ic$. At last, we denote by $\varpi$ the $4$-tuple $(I, (f_i)_{i \in I}, (s_i)_{i \in \Ic}, (j_i)_{i \in \Ic})$, and $\Cvarpi$ is a rational number (see \cite{CFRalgo}, \S 4.1). 

By convention, we let $\Cvarpi = 0$ if $\varpi = (I, (f_i)_{i \in I}, (s_i)_{i \in \Ic}, (j_i)_{i \in \Ic})$ but at least one among the $f_i$, $s_i$, $j_i$ does not lie in the above-mentioned range (for instance if $s_i > A_i$ for some $i$). This allows us to forget about the exact range of summation in \eqref{eqvarpi}.

\begin{Th} \label{thnv} 
In Theorem \ref{thzl}, for any $s_1, \ldots, s_p$ the coefficient of $\zeta(s_1, \ldots, s_p)$ in the linear combination of multiple zeta values that represents \eqref{eq7}  is equal to 
$$\sum_{ j_1 = 0} ^{n_1} \ldots \sum_{ j_p = 0} ^{n_p} \Cvidesj.$$
\end{Th}

\begin{Remarque} \label{remsun}
If $s_1 = 1$, Theorem \ref{thzl} asserts that $\zeta(s_1, \ldots, s_p)$ appears with a zero coefficient. This is consistent with Theorem \ref{thnv}, since   the assumption $\deg_{k_1}P \leq A_1(n_1+1) - 2$ yields for any $s_2, \ldots, s_p$, $j_2, \ldots, j_p$ (as in \cite{CFRsym}, \S 4.4):
$$\sum_{j_1= 0} ^{n_1} \Cvidesjsun = 0.$$
\end{Remarque}

\bigskip

\Dem of Theorem \ref{thnv}: We follow the proof \cite{ZlobinZametki2005} of Theorem \ref{thzl}. Let $z \in \C$ be such that $|z| < 1$. For any $\varpi$, let 
$$S_\varpi(z) = \sum_{k_1\geq \ldots \geq k_p \geq 1}
\frac{\prod_{i \in I} k_i^{f_i}}{\prod_{i \in \Ic}(k_i+r_i+ j_i)^{s_i}} z^{k_1},$$
so that \eqref{eq7} is the limit of $\sum_\varpi \Cvarpi S_\varpi(z)$ as $z$ tends to 1. Of course, for some $\varpi$ the function $S_\varpi(z)$ may be divergent at $z=1$, but this linear combination does have a limit thanks to \eqref{eqCV}. 

For any $\varpi$, we have an equality \cite{ZlobinZametki2005} 
$$S_\varpi(z) = \sum_{(\sigma_1, \ldots, \sigma_q)} \Domsigz \La_{\sigma_1, \ldots, \sigma_q}(z)$$
where
$\La_{\sigma_1, \ldots, \sigma_q}(z) =  \sum_{k_1\geq \ldots \geq k_q \geq 1}\frac{z^{k_1}}{k_1^{\sigma_1} \ldots k_q^{\sigma_q}}$ and $\Domsigz$ is a rational function of $z$.

The new point is that when $q=p$ we have $\Domsigz=0$  except when $I = \emptyset$ and $\sigma_1 = s_1$, \ldots, $\sigma_p = s_p$ (where $\varpi = (I, (f_i)_{i \in I}, (s_i)_{i \in \Ic}, (j_i)_{i \in \Ic})$); and in this case $\Domsigz=z^{-r_1-j_1}$.  This remark follows from the proof of \cite{ZlobinZametki2005} (see also Th\'eor\`eme 5 of \cite{CFRalgo}). 

Ending the proof as in \cite{ZlobinZametki2005}, we deduce Theorem \ref{thzl} with the additional property stated in Theorem \ref{thnv}.

\subsection{A Consequence: the Derivation Procedure} \label{subsecderiv}

Let us recall the classical ``derivation procedure'' in depth 1 (used for instance in \cite{vingtetun} and \cite{Zudilinonze}). Let $R(k) = P(k) / (k)_{n+1}^A$ be a rational fraction, with $P \in \Q[k]$ of degree at most $A(n+1)-1$, and $\ell \geq 1$. Then $\sum_{k \geq 1} R^{(\ell)}(k)$ is a linear form (over $\Q$) in $1$, $\zeta(1+\ell)$,  $\zeta(2+\ell)$, \ldots,  $\zeta(A+\ell)$. The proof of this fact is easy, by differentiating $\ell$ times the partial fraction expansion of $R(k)$.

Thanks to Theorem \ref{thnv}, this fact generalizes easily to the depth $p$ part of the linear combination that represents \eqref{eq7}. Namely, let 
$$
R(\kunp) = 
\frac{P(k_1,  \ldots, k_p)}{(k_1+r_1)_{n_1+1}^{A_1}  \ldots (k_p+r_p)_{n_p+1}^{A_p}} = 
\sum_{\varpi} 
 \Cvarpi \frac{\prod_{i \in I} k_i^{f_i}}{\prod_{i \in \Ic}(k_i+r_i+ j_i)^{s_i}},
 $$
and $\ellunp \geq 0$ be integers. Then applying $(\frac{\partial}{\partial k_1})^{\ell_1} \ldots (\frac{\partial}{\partial k_p})^{\ell_p}$ yields the partial fraction expansion of $(\frac{\partial}{\partial k_1})^{\ell_1} \ldots (\frac{\partial}{\partial k_p})^{\ell_p}R(\kunp) $, in which the polar part
$$\frac{1}{(k_1+r_1+j_1)^{s_1}  \ldots (k_p+r_p+j_p)^{s_p} }$$
appears with a non-zero coefficient only if $1 + \ell_i \leq s_i \leq A_i + \ell_i$ for any $i \in \unp$. Theorem \ref{thnv} implies that   $\zeta(\sunp)$ may appear only for these values of $\sunp$ in the linear combination of Theorem \ref{thzl} that represents
$$ \sum_{k_1\geq \ldots \geq k_p \geq 1} (\frac{\partial}{\partial k_1})^{\ell_1} \ldots (\frac{\partial}{\partial k_p})^{\ell_p}R(\kunp) , $$
assuming that it converges. However, the argument does not easily generalize to $\zeta(\sunq)$ with $q \leq p-1$. It would be interesting to investigate in this direction.

\section{Symmetry Properties of the Depth $p$ Part} \label{secprinc}

This section is the heart of the present paper. We recall the symmetry property in depth 1 connected to (very) well-poised   hypergeometric series (\S \ref{subsecdep1}), which is the origin of this work. Then we define (\S \ref{subsecenplus}) two linear representations of the group $G = \zdzp \psd \spp$: one involves rational functions in $p$ variables, and the other one formal symbols corresponding to multiple zeta values. These representations allow us to state our main result (Theorem~\ref{threp}), namely: if $R(\kunp)$ satisfies some symmetry property, then the depth $p$ part of the linear combination in Theorem \ref{thzl} satisfies a corresponding symmetry property. We prove this result in \S \ref{subsecdemprinc}, and derive several consequences of it in \S\S \ref{subsecwp} to \ref{subseccyclic}, namely special cases which yield concrete statements. 

\subsection{The Case of Depth 1} \label{subsecdep1}

Let us recall   the symmetry property used by Rivoal \cite{RivoalCRAS} and Ball-Rivoal \cite{BR} to prove that $\zeta(2n+1)$ is irrational for infinitely many integers $n$. When $e=1$, this is Assertion $(ii)$ in the Introduction.

\medskip

Let $A,n \geq 0$, $e \in \{0,1\}$, $P \in \Q[k]$ and $R(k) = P(k) / (k)_{n+1}^A$. Then the following three assertions are equivalent:
\begin{itemize}
\item $R(-k-n) = (-1)^e R(k)$,
\item $P(-k-n) = (-1)^{A(n+1)+e} P(k)$,
\item $P$ is a linear combination (over $\Q$) of $(k+\frac{n}{2})^f$ with $f \congru A(n+1) + e \bmod 2$.
\end{itemize}
If they are satisfied and $\deg P \leq A(n+1)-2$, then $\sum_{k \geq 1} R(k) $ is a linear combination, with rational coefficients, of 1 and $\zeta(s)$ with  $ 2 \leq s \leq A$ and $s \congru e \bmod 2$.

\medskip

In the next section, we generalize this symmetry property to the case of $p$ variables $\kunp$. It involves a more complicated group action, since once may permute these variables and/or make a change similar to $k \mapsto -n-k$ with respect to some of them.

\subsection{Notation and Statement of the Main Result} \label{subsecenplus}

Throughout this section, we fix integers $p \geq 1$ and $A_1, \ldots, A_p, n_1, \ldots, n_p,  r_1, \ldots, r_p \geq 0$. For any $i \in \unp$, we let 
$$K_i = k_i + r_i + \frac{n_i}{2}$$
(though $n_i$ is not assumed to be even), so that
$$(k_i+r_i)_{n_i+1} = (K_i - \frac{n_i}{2})_{n_i+1} = (K_i - \frac{n_i}{2})(K_i - \frac{n_i}{2}+1) \ldots (K_i + \frac{n_i}{2}-1)(K_i + \frac{n_i}{2})$$
is an even (resp. odd) function of $K_i$ if $n_i$ is odd (resp. even). 

\bigskip

We now consider permutations\footnote{It should be noticed that we permute the variables $K_i$; this reduces to permuting the variables $k_i$ if, and only if, $r_i + n_i/2$ is independent from $i$.}
 of the variables $\Kunp$, and changes of signs $K_i \mapsto -K_i$. We assume neither $n_1 =\ldots  = n_p$ nor $A_1 = \ldots  = A_p$, though when a permutation $K_i \mapsto K_j$ comes really into the play, the most interesting case is $n_i = n_j$ and $A_i = A_j$ (so that $(K_i - \frac{n_i}{2})_{n_i+1}^{A_i}$ maps to $(K_j - \frac{n_j}{2})_{n_j+1}^{A_j}$ and the symmetry property of the rational function $R(\kunp)$ defined in \eqref{eqRint} can be easily stated in terms of $P$). 

\bigskip

We shall denote by the same letter (e.g., $P$) a function of $\kunp$ and the corresponding function of $\Kunp$.   For instance, if $P(k_1, k_2) = (k_1 + r_1 + \frac{n_1}{2})(k_2 + r_2 + \frac{n_2}{2})$ then we let $P(K_1,K_2) = K_1K_2$.  The symmetry properties are more easily written in terms of the variables $K_i$, but we shall often translate them in terms of $k_i$. For instance, for $i \in \unp$ and $\eps \in \{-1,1\}$, the relation 
$$P(K_1, \ldots, K_{i-1}, -K_i, K_{i+1}, \ldots, K_p) = \eps P(\Kunp)$$
is equivalent to 
 $$P(k_1, \ldots, k_{i-1}, -k_i -2 r_i - n_i , k_{i+1}, \ldots, k_p) = \eps P(\kunp).$$
In the same way, when $p=2$, the assumption $P(K_2, K_1) = - P(K_1,K_2)$ is equivalent to Equation \eqref{eqPsymint}. 

\bigskip

 The permutation group  $\spp$ acts on $\zdzp$ by
 $$\gamma \cdot (\epsunp) = (\epsgammuunp)$$
 for $\gamma \in \spp$ and $(\epsunp) \in \zdzp$. This is a left action, that is 
 $\gamma \cdot (\gamma' \cdot (\epsunp)) = (\gamma \gamma') \cdot (\epsunp)$; all group actions we consider throughout this text are left actions.  This allows one to define the {\em semi-direct product} $G = \zdzp \psd \spp$ as the set-theoretic cartesian product $\zdzp \times \spp$ equipped with the law 
 $$(\epsunp, \gamma) (\epsprunp, \gamma') = (\eps_1 \eps'_{\gamma^{-1}(1)}, \ldots, \eps_p \eps'_{\gamma^{-1}(p)},  \gamma \gamma')$$
 where all group laws (including the one of $\zdz$) are written multiplicatively. A generic element of $G$ is denoted by either $(\epsunp, \gamma)$ or $(\epssoul, \gamma)$, where $\epssoul \in \zdzp$  stands for $(\epsunp)$.

\bigskip

Let $\K$ be any subfield of $\C$. We let $\VK$ be the $\K$-vector space  of all rational functions $R(\kunp)$ that can be written as 
$$R(\kunp) = \frac{P(\kunp) }{(k_1+r'_1)_{n'_1+1}^{A'_1}  \ldots (k_p+r'_p)_{n'_p+1}^{A'_p}}$$
with $P \in \K [\kunp]$ and  $A'_1, \ldots, A'_p, n'_1, \ldots, n'_p,  r'_1, \ldots, r'_p \geq 0$. 

Recalling that  $A_1, \ldots, A_p, n_1, \ldots, n_p,  r_1, \ldots, r_p \geq 0$  are fixed throughout this section, we let $\VprK$ be the subspace  of $\VK$ consisting in all rational fractions 
$$R(\kunp) = \frac{P(\kunp) }{(k_1+r_1)_{n_1+1}^{A_1}  \ldots (k_p+r_p)_{n_p+1}^{A_p}}$$
with $P \in \K [\kunp]$ such that \eqref{eqCV} holds; this assumption on the degrees of $P$ can be stated, equivalently, as 
$$\sum_{i=1}^j \deg_{k_i}R \leq -j-1 \mbox{ for any } j \in \unp.$$

We define a group homomorphism $\ro : G \rightarrow {\rm GL}(\VK)$ (that is, a $\K$-linear representation of $G$) as follows:
$$\ro(\epsunp, \gamma)(R(\Kunp)) = R(\eps_{\gamma(1)}K_{\gamma(1)}, \ldots, \eps_{\gamma(p)}K_{\gamma(p)})$$
where $\zdz$ is seen as $\{-1,1\}$ (and we keep this convention throughout this text). Let us check that $\ro$ is indeed a group homomorphism:
\begin{eqnarray*}
\ro((\epssoul, \gamma)(\epsprsoul, \gamma'))(R) 
&=& \ro(\eps_1 \eps'_{\gamma^{-1}(1)}, \ldots, \eps_p \eps'_{\gamma^{-1}(p)},  \gamma \gamma')(R)\\
&=& R(\eps_{\gamma \gamma' (1)} \eps'_{\gamma'(1)} K_{\gamma \gamma'(1)}, \ldots, 
\eps_{\gamma \gamma' (p)} \eps'_{\gamma'(p)} K_{\gamma \gamma'(p)})\\
&=& \ro(\epssoul, \gamma)(R(\eps'_{\gamma'(1)} K_{\gamma'(1)}, \ldots, \eps'_{\gamma'(p)} K_{\gamma'(p)}))\\
&=&  \ro(\epssoul, \gamma)( \ro(\epsprsoul, \gamma')(R)). 
\end{eqnarray*}

\bigskip

We define now another representation of $G$. Let $\WK$ be the $\K$-vector space generated by the formal symbols $\zetaf(\sunp)$ for positive integers $\sunp$ (recall that $p$ is fixed); these symbols are assumed to be linearly independent over $\K$, so that they make up a basis of $\WK$. Let $\WprK$ be the $\K$-vector subspace of $\WK$ generated by the symbols $\zetaf(\sunp)$ with $s_1 \geq 2$.

We have a specialization map $\varphi : \WprK \to \C$, which is $\K$-linear and maps the ``formal'' multiple zeta value $\zetaf(\sunp)$ to the usual one $\zeta(\sunp)$ (which exists since $s_1 \geq 2$). This map $\varphi$ is {\em not} injective (if $p \geq 2$), since linear relations do exist between multiple zeta values of a given depth $p$ (for instance $4 \zeta(2,4) + 13 \zeta(4,2) - 18\zeta(3,3)  = 0$). 

Let us define a representation $\tildero : G \to {\rm GL} (\WK)$ by linearity as follows:
$$\tildero(\epsunp, \gamma)(\zetaf(\sunp)) = \eps_1^{s_{\gamma^{-1}(1)}} \ldots  \eps_p^{s_{\gamma^{-1}(p)}} \zetaf(s_{\gamma^{-1}(1)}, \ldots, s_{\gamma^{-1}(p)}).$$
It might be useful to notice that $s_{\gamma^{-1}(1)}, \ldots, s_{\gamma^{-1}(p)}$ appear here (as in the action of $\spp$ on $\zdzp$) since $\zetaf(\sunp)$ behaves like a point in $\N^p$; on the contrary, $K_{\gamma(1)}, \ldots, K_{\gamma(p)}$ are used in the definition of $\ro$ since $R$ is a function on $\C^p$. 

Let us  check that $\tildero$ is  a representation:
\begin{eqnarray*}
&&\tildero((\epssoul, \gamma)(\epsprsoul, \gamma'))(\zetaf(\sunp)) \\
&=& \tildero(\eps_1 \eps'_{\gamma^{-1}(1)}, \ldots, \eps_p \eps'_{\gamma^{-1}(p)},  \gamma \gamma')(\zetaf(\sunp))\\
&=& \eps_1^{s_{{\gamma'}^{-1} \gamma^{-1}(1)}} {\eps'}_{\gamma^{-1}(1)}^{s_{{\gamma'}^{-1} \gamma^{-1}(1)}} \ldots \eps_p^{s_{{\gamma'}^{-1} \gamma^{-1}(p)}} {\eps'}_{\gamma^{-1}(p)}^{s_{{\gamma'}^{-1} \gamma^{-1}(p)}}    \zetaf(s_{{\gamma'}^{-1} \gamma^{-1}(1)} , \ldots, s_{{\gamma'}^{-1} \gamma^{-1}(p)} )\\
&=& \Big( \prod_{j=1}^p {\eps'}_j^{s_{{\gamma'}^{-1}(j)}}\Big)  \eps_1^{s_{{\gamma'}^{-1} \gamma^{-1}(1)}} \ldots  \eps_p^{s_{{\gamma'}^{-1} \gamma^{-1}(p)}}  \zetaf(s_{{\gamma'}^{-1} \gamma^{-1}(1)} , \ldots, s_{{\gamma'}^{-1} \gamma^{-1}(p)} )\\
&=&  \tildero(\epssoul, \gamma)(  {\eps'}_1^{s_{{\gamma'}^{-1} (1)}} \ldots {\eps'}_p^{s_{{\gamma'}^{-1} (p)}} \zetaf(s_{{\gamma'}^{-1}  (1)} , \ldots, s_{{\gamma'}^{-1} (p)} )\\
&=&  \tildero(\epssoul, \gamma)( \tildero(\epsprsoul, \gamma')(\zetaf(\sunp))).
\end{eqnarray*}

\bigskip

It should be noticed that $\tildero$ can not induce (via $\varphi$) a representation of $G$ on the vector space generated by the usual multiple zeta values of depth $p$, since $\tildero$ does not preserve $\Q$-linear relations between them. For instance, we have  $4 \zeta(2,4) + 13 \zeta(4,2) - 18\zeta(3,3)  = 0$ but  $4 \zeta(4,2) + 13 \zeta(2,4) - 18\zeta(3,3)  \neq 0$.

\bigskip

Recall that a {\em character} of a group $H$ is a group homomorphism from $H$ to  $\C \etoile$. We can now state our main  result, which is a complement to Theorem \ref{thzl}. 

\begin{Th} \label{threp}
Let $H$ be a subgroup of $G$, and $\chi$ be a character of $H$. In Theorem \ref{thzl}, assume that
$$R(\kunp) = \frac{P(\kunp) }{(k_1+r_1)_{n_1+1}^{A_1}  \ldots (k_p+r_p)_{n_p+1}^{A_p}}$$
satisfies
$$\ro(g) (R) = \chi(g) R $$
for any $g \in H$. Then in the linear combination that represents \eqref{eq7}, we may assume that the depth $p$ part can be written as $\varphi(x)$ for some $x \in \WprK$ such that
$$\tildero(g) (x) = \chi(g) x $$
for any $g \in H$.
\end{Th}

One may notice that the assumption on $R$ implies $\chi(H) \subset \K \etoile$, since $R \in \K(\kunp)$. 

\smallskip

\begin{Exemple} \label{exemplesym}
Let us consider the case $\K = \Q$, $n_1 = \ldots = n_p = n$, $A_1 =\ldots = A_p = A$, $H=G$ and $\chi(\epsunp, \gamma) = \eps_1 \ldots \eps_p \eps_\gamma$ where $\eps_\gamma$ is the signature of $\gamma$. Then $\ro(g) (R) = \chi(g) R $ for any $g \in H$ if, and only if, $P(\kunp)$ belongs to the set denoted by  ${\cal A}_p$ in \cite{CFRsym} (\S 2.2). The conclusion of Theorem \ref{threp} is equivalent (proceeding as in the proof of Theorems \ref{thwpprofmaxi}, \ref{thsigdel} and \ref{thcyclic} below) to the fact that the depth $p$ part of the linear combination of Theorem \ref{thzl} is a linear combination over $\Q$ of ``antisymmetric multiple zeta values'' (as defined in \cite{CFRsym})
$$\zetaantisym (s_1, \ldots, s_p) =  \sum_{\gamma\in \spp} \eps_{\gamma} 
\zeta(s_{\gamma(1)}, \ldots, s_{\gamma(p)}) ,$$
in which $\sunp \geq 3$ are odd.  Actually this property follows from Th\'eor\`eme 4 of \cite{CFRsym}, which provides much more information. It would be very interesting, for other pairs $(H,\chi)$, to generalize Theorem \ref{threp} and obtain some information about the whole linear combination of Theorem \ref{thzl}, and not only its depth $p$ part. This would be specially interesting from the diophantine point of view, since the main drawback of Th\'eor\`eme 4 of \cite{CFRsym} is that $H$ is too big (therefore a lot of constraints have to be imposed on $P$). 
\end{Exemple}

\subsection{Proof of the Main Result} \label{subsecdemprinc}

In this paragraph, we prove Theorem \ref{threp} by connecting the two group actions defined in \S \ref{subsecenplus} thanks to an equivariant linear map $f$ (see Proposition \ref{prop30} below). Let us define $f : \VK \to \WK$ now. Let $R \in \VK$, and $( \Cvarpi )$ be the family of coefficients defined by
\begin{equation} \label{eqvarpide}
R(k_1,  \ldots, k_p) = 
\sum_{\varpi} 
 \Cvarpi \frac{\prod_{i \in I} k_i^{f_i}}{\prod_{i \in \Ic}(k_i+r_i+ j_i)^{s_i}},
\end{equation}
with $\varpi = (I, (f_i)_{i \in I}, (s_i)_{i \in \Ic}, (j_i)_{i \in \Ic})$ as in Equation \eqref{eqvarpi}.
Then we let 
$$f(R) = \sum_{\sunp} \Big( \sum_{ j_1 = 0} ^{n_1} \ldots \sum_{ j_p = 0} ^{n_p} \Cvidesj \Big) \zetaf(\sunp)$$
where in the sum, $s_i$ ranges from $1$ to $A_i$; but this range is not very important, since $\Cvidesj$ is zero otherwise. 

Remark \ref{remsun} proves that  $f(R) \in \WprK$ as soon as $R \in \VprK$ (which is the case in Theorem \ref{threp}). When $R \in \VprK$, Theorem \ref{thnv} asserts that $\varphi(f(R))$ is the depth $p$ part of the linear combination constructed in  Theorem \ref{thzl}. Now the key point in the proof of Theorem \ref{threp} is the following result:

\begin{Prop} \label{prop30} The $\K$-linear map $f  : \VK \to \WK$ is equivariant with respect to the actions of $G$, that is
$$f \circ \ro(g) = \tildero(g) \circ f$$
holds for any $g \in G$.
\end{Prop}

To deduce Theorem \ref{threp} from this proposition, let $ R \in \VprK$ be such that $\ro(g) (R) = \chi(g) R $ for any $g \in H$. Since $\chi(H) \subset \K \etoile$ (as noticed after the statement of Theorem \ref{threp}), Proposition \ref{prop30} yields $ \tildero(g)(f(R)) = \chi(g) f(R)$ for any $g \in H$, by $\K$-linearity. Letting $x = f(R) \in \WprK$ concludes the proof of  Theorem \ref{threp}.

\bigskip

\Dem of  Proposition \ref{prop30}: Let $g = (\epsunp, \gamma) \in G$, $R \in \VprK$ and consider $R' = \ro(g) (R) $. Denote by $( \Cprvarpi )$  the family of coefficients associated with $R'$ as in Equation \eqref{eqvarpide}. By definition of $\ro$ we have
$$ R'(\Kunp) = R(\eps_{\gamma(1)}K_{\gamma(1)}, \ldots, \eps_{\gamma(p)}K_{\gamma(p)}).$$
Expanding both sides into partial fractions as in  \eqref{eqvarpide} yields
\begin{equation} \label{eqpa}
 \sum_{\varpi} 
 \Cprvarpi \frac{\prod_{i \in I} (K_i - r_i - \frac{n_i}{2})^{f_i}}{\prod_{i \in \Ic}(K_i - \frac{n_i}{2}+ j_i)^{s_i}}
=
\sum_{\varpi} 
 \Cvarpi \frac{\prod_{i \in I} (\eps_{\gamma(i)} K_{\gamma(i)} - r_i -  \frac{n_i}{2})^{f_i}}{\prod_{i \in \Ic}(\eps_{\gamma(i)} K_{\gamma(i)} - \frac{n_i}{2}+ j_i)^{s_i}}.
\end{equation}
Now we let $\zdz = \{-1,1\}$ act on $\zeron$ by 
$$
\left\{
\begin{array}{l}
\eps \cdot j =j \mbox{ for } \eps =1 , \\
\eps \cdot j =n-j \mbox{ for } \eps = -1,
\end{array}
\right.
$$
in such a way that 
\begin{equation}\label{eqeps}
( -\frac{n}{2} + j) \eps =  -\frac{n}{2} +\eps \cdot  j
\end{equation}
for any $\eps \in \{-1,1\}$ and $j \in \zeron$. Of course this action depends on $n$, but it is clear from the context that $\eps \cdot j_i$ refers to the case $n= n_i$ (even if $\eps$ is denoted by $\eps_\ell$ for some $\ell$). 

Now the uniqueness of the partial fraction expansion \eqref{eqpa} yields, for any $\ssoul = ( \sunp)$ and  $\jsoul = (\junp)$ (using \eqref{eqeps}):
$$
 \frac{\Cprvidesjsoul}{\prod_{i= 1}^p  (K_i - \frac{n_i}{2} + j_i)^{s_i}}
=
\frac{\Cvidesjgam}{\prod_{\ell = 1}^p \eps_\ell ^{s_\ell} (K_\ell - \frac{n_\ell}{2} + j_\ell)^{s_\ell}},
$$
that is
\begin{equation}\label{eqC}
 \Cprvidesjsoul
 = \eps_1^{s_1} \ldots  \eps_p^{s_p} \Cvidesjgam .
\end{equation}
Now the equality
$$f(\ro(g)(R)) = f(R') =   \sum_{\sunp \geq 1} \Big( \sum_{ j_1, \ldots, j_p}   \Cprvidesj \Big) \zetaf(\sunp)$$
yields, thanks to \eqref{eqC}:
\begin{eqnarray*}
f(\ro(g)(R)) &=& \sum_{\sunp \geq 1} \Big( \sum_{ j_1, \ldots, j_p}  \Cvidesjgam \Big)\\
&& \quad \quad \quad \quad \quad \quad \quad \quad  \quad \quad \quad \quad   \eps_1^{s_1} \ldots  \eps_p^{s_p}  \zetaf(\sunp).
\end{eqnarray*}
Letting $s'_i = s_{\gamma(i)}$ and $j'_i = \eps_i \cdot j_{\gamma(i)}$  for any $i \in \{1, \ldots, p\}$, this equality reads
\begin{eqnarray*}
f(\ro(g)(R)) &=& \sum_{s'_1, \ldots, s'_p \geq 1} \Big( \sum_{ j'_1, \ldots, j'_p} \Cvidesprjpr \Big)\\
&& \quad \quad \quad \quad \quad \quad \quad \quad  \quad \quad \quad \quad   \eps_1^{s'_{\gamma^{-1}(1)}} \ldots  \eps_p^{s'_{\gamma^{-1}(p)}}   \zetaf(s'_{\gamma^{-1}(1)}, \ldots, s'_{\gamma^{-1}(p)}).
\end{eqnarray*}
By definition of $\tildero $ this means
$f(\ro(g)(R)) =  \tildero(g)( f(R))$, thereby concluding the proof of Proposition \ref{prop30}.

\subsection{A Consequence Involving the Parity of $s_i$} \label{subsecwp}

Let us start with the following consequence of Theorem \ref{threp}, in which no permutation of the variables $\Kunp$ is involved. 

\begin{Th} \label{thwpprofmaxi}
In the situation of Theorem \ref{thzl}, assume that for some integers
 $e_1, \ldots, e_p$ we have
 $$P(k_1, \ldots, k_{i-1}, -k_i -2 r_i - n_i , k_{i+1}, \ldots, k_p) = (-1)^{A_i(n_i+1)+e_i} P(\kunp)$$
 for any $i \in \unp$. Then in the linear combination that represents \eqref{eq7}, any multiple zeta value $\zeta(\sunp)$ of depth $p$ that appears with a non-zero coefficient satisfies
 $$s_i \congru e_i \mod 2 \mbox{ for any } i \in \unp.$$
\end{Th}

This result generalizes the symmetry phenomenon used by Rivoal \cite{RivoalCRAS} and Ball-Rivoal \cite{BR} to prove that $\zeta(2n+1)$ is irrational for infinitely many integers $n$; namely this property (recalled in \S \ref{subsecdep1}) is obtained for $p=1$. However, Theorem \ref{thwpprofmaxi} is not as powerful as  the results of  \cite{CFRsym}, since it concerns only the depth $p$ part of the linear combination (a good challenge is to strengthen it: see \S \ref{secconj}). 

\bigskip

For instance, when $p=2$ and $r_1 \geq r_2 + n_2  + 1$,  Theorem \ref{thwpprofmaxi} yields linear forms in 1, $\zeta(s)$ with $2 \leq s \leq \max(A_1, A_2)$, and $\zeta(s_1,s_2)$ with $1 \leq s_i \leq A_i$ ($i \in \{1,2\}$) and $s_i$ of fixed parity. If in addition $A_1 = A_2 = 3$   and $\deg_{k_i}P \leq 3n_i + 1$, one obtains $1$, $\zeta(2)$, $\zeta(3)$, and exactly one  multiple zeta value among $\zeta(2,2)$, $\zeta(2,3)$, $\zeta(3,2)$, $\zeta(3,3)$. More generally, plugging this symmetry phenomenon into Theorem \ref{th1} enables one to get only one multiple zeta value of weight $p$: this is Theorem \ref{th2} stated in the Introduction. It would be very interesting to obtain analogous symmetry properties on $P$ that ensure that only some multiple zeta values of weights $< p$ appear; but this seems to be a difficult question. Some conjectures in this direction are made in the last section of this paper.

\bigskip

\Dem of Theorem \ref{thwpprofmaxi}: Let $H = \zdzp \times \{ \Id\}$ and $\chi(\epsunp, \Id) = \eps_1^{e_1} \ldots \eps_p ^{e_p}$. The assumption on $P$ means 
$\ro(g) (R) = \chi(g)R$ for any $g \in H$, so that Theorem \ref{threp} applies. To conclude the proof, what remains is to understand which elements $x \in \WprK$ satisfy 
$\tildero(g) (x) = \chi(g)x$ for any $g \in H$. Writing 
$$x = \sum_{\sunp} \lambda[\sunp] \zetaf(\sunp)$$
with $\lambda[\sunp] \in \K$, this condition means 
$$\eps_1^{s_1} \ldots \eps_p ^{s_p} \lambda[\sunp]   =  \eps_1^{e_1} \ldots \eps_p ^{e_p} \lambda[\sunp] 
$$
for any $\epsunp$. This is equivalent to 
$$\lambda[\sunp] = 0 \mbox{ as soon as } s_i \not\congru e_i \mod 2 \mbox{ for some } i \in \unp,$$
thereby concluding the proof of Theorem \ref{thwpprofmaxi}. 

\bigskip

\begin{Remarque} If the symmetry assumption on $P$ in Theorem \ref{thwpprofmaxi} is satisfied only for some values of $i \in \unp$, then the conclusion on the parity of $s_i$ holds for these values of $i$. This can be proved in the same way, or deduced from Theorem \ref{thwpprofmaxi} by decomposing $P(K_1, \ldots, K_p)$ into even and odd parts with respect to the variables $K_i$ for which no symmetry is assumed.
\end{Remarque}

\subsection{Other Results in Depth Two} \label{subsecdep2}

In this section, we restrict to the case $p=2$ and allow the variables $K_1$, $K_2$ to be permuted. To make the statements simpler, we consider only monomials $P$ in the numerator. However, the possible diophantine applications would come from suitable linear combinations of these monomials (with fixed parity conditions on $e$ and $f$). 

\begin{Th} \label{thsigdel} Let $A , n, r, e, f$ be non-negative integers such that 
$  e+f \leq A(n+1)-2 $.

Let $S = A$ if $r \geq n+1$, and $S  = 2A$ otherwise.
Then the series
\begin{equation} \label{eqsigdel}
\sum_{k_1 \geq k_2 \geq 1} 
\frac{(k_1 + k_2 + n + r)^e (k_1 - k_2  + r)^f}{(k_1+r)_{n+1}^A (k_2)_{n+1}^A}
\end{equation}
is a linear combination with rational coefficients  of 1, $\zeta(s)$ with $2 \leq s \leq S$, and:
\begin{itemize}
\item if $e$ and $f$ are even,   $\zeta(s,s') + \zeta(s',s)$ with $2 \leq s \leq s' \leq A$ and $s' \congru s \bmod 2$.
\item if $e$ is even and $f$ is odd,   $\zeta(s,s') - \zeta(s',s)$ with $2 \leq s < s' \leq A$ and $s' \not\congru s \bmod 2$.
\item if $e$ is odd and $f$ is even,   $\zeta(s,s') + \zeta(s',s)$ with $2 \leq s < s' \leq A$ and $s' \not\congru s \bmod 2$.
\item if $e$ and $f$ are odd, $\zeta(s,s') - \zeta(s',s)$ with $2 \leq s < s' \leq A$ and $s' \congru s \bmod 2$.
\end{itemize}
\end{Th}

The identity $\zeta(s,s') + \zeta(s',s) = \zeta(s) \zeta(s') - \zeta(s+s')$ may be used (when $f$ is even) to express in a different way the conclusion of this theorem. We obtain for instance the following corollary:

\begin{Cor}
Under the assumptions of Theorem \ref{thsigdel},  if $e$ is odd,  $f$ is even and $r \geq n+1$, then \eqref{eqsigdel} is  a linear combination over the rationals of:
\begin{itemize}
\item 1, 
\item $\zeta(s)$ with $2 \leq s \leq A$, 
\item $\zeta(s)$ with $A+1 \leq s \leq 2A-1$ and $s$ odd, 
\item $\zeta(s) \zeta(s')$ with $2 \leq s < s' \leq A$ and $s' \not\congru s \bmod 2$.
\end{itemize}
\end{Cor}

If $A=3$, this corollary yields linear forms in $1$, $\zeta(2)$, $\zeta(3)$ and $\zeta(2) \zeta(3)  - \zeta(5)$. 

\bigskip

When $f$ is odd, Theorem \ref{thsigdel} yields antisymmetric multiple zeta values of depth 2 (as defined in \cite{CFRsym}), that is $\zeta(s,s') - \zeta(s',s)$. According to the parity of $e$, we know whether $s$ and $s'$ have the same parity or not. But even when $e$ is odd (so that $s' \congru s \bmod 2$), $s$ may be even or odd (however, see Corollary \ref{corcommun2} below). This is an important difference with Th\'eor\`eme 3 of \cite{CFRsym}, where $\zeta(s,s') - \zeta(s',s)$ appears only when $s$ and $s'$ are odd. Another difference is that $\zeta(s)$ appears in Theorem \ref{thsigdel} for any $s \leq S$, whereas it does in Th\'eor\`eme 3 of \cite{CFRsym} for odd values of $s \leq 2A-1$.

\bigskip

\Dem of Theorem \ref{thsigdel}:   Let $\tau \in \sdeux$ be the transposition, and $H$ be the subgroup of $\zdzd \psd \sdeux$ generated by $(1,1,\tau)$ and $(-1,-1,\Id)$. Let $\chi : H \to \{-1,1\}$ be the character defined by $\chi(1,1,\tau) = (-1)^f$ and $\chi(-1,-1,\Id) = (-1)^{e+f}$. The rational function
$$R(K_1, K_2) = \frac{(K_1+K_2)^e (K_1 - K_2)^f}{(K_1 - \frac{n}{2})_{n+1}^A (K_2 - \frac{n}{2})_{n+1}^A }$$
satisfies the symmetry properties
\begin{equation} \label{eqsymsigdel}
R(K_2, K_1)  = (-1)^f R(K_1,K_2) \mbox{ and } R(-K_1, -K_2) = (-1)^{e+f} R(K_1,K_2).
\end{equation}
This means exactly that $\ro(g)(R)  = \chi(g) R$ holds for the two above-mentioned generators of $H$; therefore this relation holds for any $g \in H$, and Theorem \ref{threp} applies (with $\K = \Q$). 

Let
$$x = \sum_{s_1, s_2} \lambda[s_1, s_2] \zetaf(s_1, s_2)$$
be an element of $\WprQ$:  the rational number  $\lambda[s_1, s_2] $ vanishes as soon as $s_1=1$. 
The relations $\tildero(1,1,\tau)(x) = (-1)^f x$ and $\tildero(-1,-1,\Id)(x) = (-1)^{e+f} x$ mean, respectively, $\lambda[s_1, s_2] = (-1)^f \lambda[s_2, s_1]$ and
$\lambda[s_1, s_2] ((-1)^{e+f} - (-1)^{s_1+s_2}) = 0$ for any $s_1, s_2$. They imply $\lambda[s_1, s_2] = 0$ when $s_1 + s_2 \not\congru e+f \bmod 2$, so that $x$ is a linear combination of $\zetaf(s_1,s_2) + (-1)^f \zetaf(s_2,s_1)$ with $s_1 \congru s_2 + e+f \bmod 2$ and $s_1, s_2 \geq 2$. This concludes the proof of Theorem \ref{thsigdel}. 

\medskip

\begin{Remarque} In the statement of Theorem \ref{thsigdel}, one could have replaced the specific form of $P(k_1, k_2)$ by the assumption \eqref{eqsymsigdel}. 
\end{Remarque}

\bigskip

When $e=f$ in \eqref{eqsigdel}, one may apply either Theorem \ref{thsigdel} or Theorem \ref{thwpprofmaxi}. Since the linear combination in the conclusion of both is the same (namely the one constructed in the proof of Theorem \ref{thzl}, see Remark \ref{remunicite}), we derive the following corollary.

\begin{Cor} \label{corcommun2}
Let $A , n, r, e $ be non-negative integers such that 
$  2e \leq A(n+1)-2 $.

Let $S = A$ if $r \geq n+1$, and $S  = 2A$ otherwise.
Then the series
\begin{equation} \label{eqcommun2}
\sum_{k_1 \geq k_2 \geq 1} 
\frac{\Big((k_1 + k_2 + n + r) (k_1 - k_2  + r)\Big)^e}{(k_1+r)_{n+1}^A (k_2)_{n+1}^A}
\end{equation}
is a linear combination with rational coefficients of 1, $\zeta(s)$ with $2 \leq s \leq S$, and:
\begin{itemize}
\item if $e$ is  even,   $\zeta(s,s') + \zeta(s',s)$ with $2 \leq s \leq s' \leq A$ and $s' \congru s \congru A(n+1) \bmod 2$.
\item if $e$ is  odd, $\zeta(s,s') - \zeta(s',s)$ with $2 \leq s < s' \leq A$ and $s' \congru s \congru A(n+1)  \bmod 2$.
\end{itemize}
\end{Cor}

This corollary could be deduced directly from Theorem \ref{threp}, by considering the subgroup $H$ generated by $\zdzd \times \{\Id\}$ and $(1,1,\tau)$. 

\medskip

From a diophantine point of view, this corollary seems to be more interesting when $A(n+1)$ is odd. In this case, when $e$ is odd, we get a linear form in 1, $\zeta(s)$ with $2 \leq s \leq S$, and $\zeta(s,s') - \zeta(s',s)$ with $2 \leq s < s' \leq A$ and $s$, $s'$ odd. This looks like  Th\'eor\`eme 3 of \cite{CFRsym}, with a major difference: in Corollary \ref{corcommun2}, $\zeta(s)$ appears   (in general) for both even and odd values of $s$,  whereas in Th\'eor\`eme 3 of \cite{CFRsym} only odd values of $s$ are involved. For instance, when $A=3$, $n=1$, $r=2$  and $e=1$, the series \eqref{eqcommun2} is equal to $- \frac{43}{16} + 2 \zeta(2) - \frac12 \zeta(3)  = 0.001339\ldots$. 

\begin{Remarque} \label{rem35} 
When $A \leq 3$ and $e$ is odd, Corollary \ref{corcommun2} proves that the double sum \eqref{eqcommun2} is a linear form in $1$, $\zeta(2)$, $\zeta(3)$, \ldots, $\zeta(S)$: no multiple zeta value of depth 2 appears. 
\end{Remarque}

\subsection{A Property Involving Cyclic Permutations} \label{subseccyclic}

This section is the only one where a  field $\K$ other than $\Q$ is used.

\begin{Th} \label{thcyclic} 
Let $p \geq 1$ and $A , n, \runp \geq 0$ be   integers, and $\xi \in \C\etoile$ be  such that 
$\xi^p = 1$. Let $\K = \Q(\xi)$, and $P \in \K [\kunp]$ be such that \eqref{eqCV} holds and 
$$P(K_2,K_3,  \ldots, K_p, K_1) = \xi \,  P(\Kunp).$$
Then the series
$$
\sum_{k_1\geq \ldots \geq k_p \geq 1}
\frac{P(k_1,  \ldots, k_p)}{(k_1+r_1)_{n+1}^{A}  \ldots (k_p+r_p)_{n+1}^{A}}
$$
is a linear combination over $\K$ as in Theorem \ref{thzl}, of which the depth $p$ part is a linear combination over $\K$ of 
$$\sum_{i=1}^p \xi^{i-1} \zeta(s_i, s_{i+1}, \ldots, s_p, s_1, s_2, \ldots, s_{i-1})$$
with $\sunp \geq 2$.
\end{Th}

This theorem can be used with $\xi = 1$, and also if $p$ is even with $\xi = -1$; in both cases $\K = \Q$. When $p=2$, $A=3$,  and $\xi = -1$, it reduces to Theorem \ref{th3} stated in the Introduction.

\bigskip

\Dem of Theorem \ref{thcyclic}: 
Let $\gamma_0 \in \spp$ be the cyclic permutation that maps 1 to 2, 2 to 3, \ldots, $p$ to 1. Denote by $H$ the subgroup of $G$ generated by $(1, \ldots, 1, \gamma_0)$; obviously $H$ is cyclic of order $p$. The assumption on $P$ means $\ro(1,\ldots, 1, \gamma_0)(R) = \xi R$ since $n_1 = \ldots  = n_p = n$ and $A_1 = \ldots = A_p = A$. This implies  $\ro(g) (R) = \chi(g)R$ for any $g \in H$,  where $\chi : H \to \K \etoile $ is the (unique) character defined by $\chi(1,\ldots, 1, \gamma_0) = \xi$. 

Therefore Theorem \ref{threp} applies; let us translate its conclusion. Let
$$x = \sum_{\sunp} \lambda[\sunp] \zetaf(\sunp)$$
be an element of $\WprK$, with $ \lambda[\sunp] \in \K$ and $\lambda[\sunp] = 0$ as soon as $s_1=1$. Then
$$\tildero(1, \ldots, 1, \gamma_0)(x) =  \sum_{\sunp} \lambda[\sunp] \zetaf(s_p, s_1, s_2, \ldots, s_{p-1}).$$
Therefore the equality $\tildero(g)(x) = \chi(g)x$, for $g = (1, \ldots, 1, \gamma_0)$, means
$$\lambda[s_2, s_3, \ldots, s_{p}, s_1] = \xi \,  \lambda[\sunp] $$
for any $\sunp$. This implies, for any $\sunp$:
\begin{eqnarray*}
&& \sum_{i=1} ^p \lambda[s_i, s_{i+1}, \ldots, s_p, s_1, s_2, \ldots, s_{i-1}] \zetaf(s_i, s_{i+1}, \ldots, s_p, s_1, s_2, \ldots, s_{i-1}) \\
& =& \lambda[\sunp]  \sum_{i=1}^p \xi^{i-1} \zeta(s_i, s_{i+1}, \ldots, s_p, s_1, s_2, \ldots, s_{i-1}),
\end{eqnarray*}
the coefficient $ \lambda[\sunp]$ being zero if $s_i = 1$ for at least one $i$. This concludes
 the proof of Theorem \ref{thcyclic}. 

\bigskip

\begin{Exemple} With $p=3$ and $\xi = (-1 + i \sqrt{3})/2$, this theorem provides linear combinations over $\Q(\xi)$ of which the depth 3 part involves
$$\zeta(s_1, s_2, s_3) + \xi \,  \zeta(s_2, s_3, s_1) + \xi^2 \, \zeta(s_3, s_1, s_2)$$
with $s_1, s_2, s_3 \geq 2$.
\end{Exemple}

\section{Conjectures on the Depth $p-1$ Part} \label{secconj}

An interesting generalization of the results proved in \S \ref{secprinc} would be to describe the part in depth $\leq p-1$ of the linear combination given by Theorem \ref{thzl}, under suitable symmetry properties of the rational function $R(\kunp)$. In the special case where $H=G$ and $\chi(\epsunp, \gamma) = \eps_1 \ldots \eps_p \eps_\gamma$, this was done in \cite{CFRsym} (see Example \ref{exemplesym} in \S \ref{subsecenplus}). It could be useful to obtain such a statement for other pairs $(H,\chi)$; in particular, if $H$ is smaller then the conditions to be imposed on $R$ are weaker, so that it is more reasonable to hope for diophantine applications. 

Maybe such a result could be obtained considering {\em coloured} multiple zeta values, and series
$$
\sum_{k_1\geq \ldots \geq k_p \geq 1}
\frac{P(k_1,  \ldots, k_p)}{(k_1+r_1)_{n_1+1}^{A_1} \ldots (k_p+r_p)_{n_p+1}^{A_p}} 
\xi_1^{k_1} \ldots \xi_p ^{k_p}
$$
where $\xiunp$ are roots of unity. The algorithm described in \cite{CFRalgo} could enable one (theoretically) to compute this kind of sums. However, the implementation \cite{CFRweb} has been done only for $\xi_1 = \ldots = \xi_p = 1$, so that we could not proceed to any experiment in the general setting.

In the situation of Theorem \ref{thzl} with $p=3$, we have found some examples of pairs $(H,\chi)$ for which the linear combination involves neither $\zeta(2,2)$ nor $\zeta(2,2,2)$. This is shown by the following theorem, that we have proved by computing all series \eqref{eqhoffman} thanks to \cite{CFRweb}. 

\begin{Th} \label{thhoffman}
Assume $A = 3$ and $n \in \{0, 1,2\}$. We let 
$$K_1 = k_1 + \frac{5n}{2} + 2, \hspace{0.3cm}
K_2 = k_2 + \frac{3n}{2}+1,  \hspace{0.3cm} 
\mbox{ and }  \hspace{0.3cm}
K_3 = k_3 + \frac{n}{2}$$
and consider a polynomial $P(K_1,K_2,K_3)$ among one of the following four families (where the exponents are non-negative integers):
\begin{enumerate}
\item \label{famiwp} $K_1^e K_2^f K_3^g $ with $e \congru g \congru A(n+1) + 1 \bmod 2$,
\item \label{famirsym} $K_2^e (K_1-K_3)^f (K_1+K_3)^g (K_1K_3)^h$ with $e \congru A(n+1) + 1 \bmod 2$ and $g \congru A(n+1) \bmod 2$, 
\item \label{famic} $K_1^e (K_2^2-K_3^2)^f (K_2K_3)^g$  with $f \congru g \congru A(n+1) + 1 \bmod 2$,
\item \label{famid} $K_3^e (K_1^2-K_2^2)^f (K_1K_2)^g$  with $f \congru g \congru A(n+1) + 1 \bmod 2$.
\end{enumerate}
If the series
\begin{equation} \label{eqhoffman}
\sum_{k_1 \geq k_2 \geq k_3 \geq 1} 
\frac{ P(K_1, K_2, K_3) }{(k_1+2n+2)_{n+1}^A (k_2+n+1)_{n+1}^A (k_3)_{n+1}^A}
\end{equation}
is convergent, then it is a linear combination (over the rationals) of multiple zeta values $\zeta(s_1, \ldots, s_q)$ with $q \in \{0, \ldots, 3\}$, $1 \leq s_i \leq 3$ for any $i$, $s_1 \geq 2$, in which neither $\zeta(2,2)$ nor $\zeta(2,2,2)$ appears.

If, in addition, \eqref{eq8} holds then it is a linear combination with rational coefficients of 1, $\zeta(2)$, $\zeta(3)$, $\zeta(2,3)$, $\zeta(3,2) $, $\zeta(3,3) $, $\zeta(2,2,3)$, $\zeta(2,3,2)$, $\zeta(2, 3,3)$,  $\zeta(3,2,2)$, $\zeta(3,2,3)$, $\zeta(3,3, 2)$, $\zeta(3,3,3)$. 
\end{Th}

The computations were too heavy for $n=3$, but we propose the following conjecture.

\begin{Conj} \label{conjhoffman}
Theorem \ref{thhoffman} holds for $A=3$ and any integer  $n \geq 0$.
\end{Conj}

\bigskip

As in Theorem \ref{thsigdel}, any diophantine application of this conjecture would involve suitable linear combinations of the peculiar polynomials $P$ defining each of the 4 families in Theorem \ref{thhoffman}.

For each of these four families,  there is a pair $(H,\chi)$ such that Theorem \ref{threp} applies, and gives some properties of the depth 3 part of the linear combination (which imply that $\zeta(2,2,2)$ does not appear). Since Theorem \ref{thhoffman} does not say more than Theorem \ref{thzl} about the depth 1 part, the open problem in Conjecture \ref{conjhoffman} is to prove that $\zeta(2,2)$ does not appear.

\bigskip

It is likely that Theorem  \ref{thhoffman} (and Conjecture \ref{conjhoffman}) can be generalized to other values of $r_1, r_2, r_3$ satisfying \eqref{eq6}, but we did not try to prove it. Another generalization would be the case $ A \geq 4$.

\bigskip

It would be interesting to obtain an analog of  Theorem  \ref{thhoffman} with $p=4$, in which $\zeta(2,2)$,  $\zeta(2,2,2)$ and $\zeta(2,2,2,2)$ disappear. It should be noted that we did not succeed in obtaining an analogous statement, with $p \in \{2,3\}$, in which $\zeta(2)$ disappears.

\bigskip

At last, another direction would be to study, in Theorem \ref{thzl}, the maximal weight part. We have found many examples in which it vanishes (as in depth 1 when $R(-n-k) = -R(k)$ and $A$ is even, see \S \ref{subsecdep1}). It is likely that a general statement can be proved, in the same spirit as the present text.

\newcommand{\url}{\texttt}

\providecommand{\bysame}{\leavevmode ---\ }
\providecommand{\og}{``}
\providecommand{\fg}{''}
\providecommand{\smfandname}{\&}
\providecommand{\smfedsname}{\'eds.}
\providecommand{\smfedname}{\'ed.}
\providecommand{\smfmastersthesisname}{M\'emoire}
\providecommand{\smfphdthesisname}{Th\`ese}

\end{document}